\documentclass{article}[12pt]
\usepackage{amsthm,amsfonts,amsmath,amssymb}
\usepackage{amsmath}
\usepackage{amsthm}
\usepackage{amsfonts}
\usepackage{amssymb}

\numberwithin{equation}{section}

\newtheorem{thm}{Theorem}[section]
\newtheorem{assertion}[thm]{Proposition}
\newtheorem{lemma}[thm]{Lemma}

\theoremstyle{remark}

\newcommand{\mmp}{\mathbb{P}}
\newcommand{\me}{\mathbb{E}}
\newcommand{\mn}{\mathbb{N}}

\def\1{{\bf 1}}

\begin{document}

\title{Small parts in the Bernoulli sieve\footnote{The work of A. Iksanov and U. Roesler was supported
by the German Scientific Foundation (project no. 436UKR
113/93/0-1).}}
\date{\today}
\author{Alexander Gnedin\footnote{e-mail address: gnedin@math.uu.nl}\\
\small{\emph{Department of Mathematics}, \emph{Utrecht University}},\\
\small{\emph{Postbus 80010, 3508 TA Utrecht, The Netherlands}}\\
Alex Iksanov\footnote{e-mail address: iksan@unicyb.kiev.ua}\\
\small{\emph{Faculty of Cybernetics},
\emph{National T. Shevchenko University}},\\
\small{\emph{01033 Kiev, Ukraine}}\ \\
Uwe Roesler\footnote{e-mail address: roesler@computerlabor.math.uni-kiel.de}\\
\small{\emph{Mathematisches
Seminar}, \emph{Christian-Albrechts-Universit\"{a}t zu Kiel}},\\
\small{\emph{Ludewig-Meyn-Str.4, D-24098 Kiel, Germany}}}

\maketitle

\begin{abstract}
\noindent
Sampling from a random discrete distribution induced by a `stick-breaking' process is considered. 
Under a  moment condition, it is shown that 
the asymptotics of the sequence of occupancy numbers,
 and of the small-parts counts (singletons, doubletons, etc) can be read off from a limiting model involving 
a unit Poisson point process and a self-similar renewal process on the halfline.
\end{abstract}

\section{Introduction}

A multiplicative renewal process (also known as residual allocation model or stick-breaking)
is a random sequence
$B=(P_j: j=0,1,\ldots )$ of the form
\begin{equation}\label{StBr}
 \ \ P_j=\prod_{i=1}^j W_i, \
\end{equation}
(so $P_0=1$) where $(W_i: i=1,2,\ldots)$ are independent copies of
a random variable $W$ taking values in $]0,1[$. We shall assume
that the support of the distribution of $W$ is not a geometric
sequence or, equivalently, that the distribution of the variable
$|\log W|$ is non-lattice, and also assume that
\begin{equation}\label{mu}
\mu:={\mathbb E}[-\log W]<\infty.
\end{equation}
The `stick-breaking' set $B$ will be viewed
as a simple point process, with $0$ being the only accumulation point.
The complement $B^c=\,[0,1]\,\setminus B$ is an open set comprised of the component intervals $]P_{j+1},P_j[$ for $j=0,1,\ldots$.

\par Let $U_1,U_2,\ldots$ be independent uniform $[0,1]$ random points, also independent of $B$,
and for each $n$ let  $U_{n,1}<\ldots<U_{n,n}$
be the order statistics of $U_1,\ldots,U_n$.
These data define a random occupancy scheme, in which a collection of $n$ `balls' $U_1,\ldots,U_n$
is sequentially sorted  into `boxes' $]P_{j},P_{j-1}[$\,, $j=1,2,\ldots$.
In the most studied and analytically best tractable case the law of $W$ is beta$(\theta,1)$, and the allocation
of `balls-in-boxes' belongs to the circle of questions around the Ewens sampling formula \cite{ABT, TSF}.
Let $K_n$ be the number of occupied `boxes' and $K_{n,r}$ be the number 
of `boxes' occupied by exactly $r$ `balls', so $\sum_{r>0}K_{n,r}=K_n$ and $\sum_{r>0} r K_{n,r}=n$.
We also define $K_{n,0}$ to be the number of unoccupied interval components of $B^c\cap \,[U_{n,1},1]$, so that
$K_n=I_n-K_{n,0}$, with
$I_n:=\min\{i:
P_i<U_{n,1}\}$ being
the index of the leftmost occupied interval.

\par  In 
\cite{Gne, GINR} the 
 renewal theory was applied to explore 
the spectrum of possible limit laws 
for $K_n$, including 
normal, stable and Mittag-Leffler distributions.
In the present note we
 focus on
the variables
$K_{n,r}$. 
We approach the $K_{n,r}$'s via the
  {\it occupancy counts}
$$Z_n^{(i)}:=\#\{1\leq j\leq n: ~U_j\in~
]P_{I_n-i+1}, P_{I_n-i}[\,\}, \ \ i=1,2,\ldots$$ in the 
left-to-right order of intervals, where we adopt the convention
$Z_n^{(i)}=0$ for $i> I_n$.
Extending a result from \cite{GINR} about $Z_n^{(1)}$, we will show that the $Z_n^{(i)}$'s
jointly converge to the sequence of occupancy numbers in a limiting model that involves 
a Poisson process and another self-similar point process on the halfline.

\par From a  viewpoint,
$B$ is an exponential transform of the range of a 
subordinator with finite L{\'e}vy measure, that is of 
a compound Poisson process.
Asymptotics of $K_n, K_{n,r}$'s have been studied in a similar occupancy model
for subordinators with infinite  L{\'e}vy measures
\cite{BG, GPY1, GPY2}.
In the infinite measure case
neither the counts $Z_n^{(i)}$ nor $I_n$ can be defined, because $B$ is then a random Cantor set, hence
there are infinitely many unoccupied intervals between any two
occupied components of $B^c$.

\section{Occupancy counts}

For $1\leq m\leq n$ the probability that the interval $]P_1,P_0[$ contains $m$ out of $n$ uniform points is
$$p(n:m)={n\choose m} \me[W^{n-m}(1-W)^m].$$
Let $(n_1,\ldots,n_k)$ be a {\it weak composition} of $n$, meaning that $n_1>0, n_2\geq 0,\ldots,n_k\geq 0$ and $n_1+\ldots+n_k=n$.
The structure (\ref{StBr}) and elementary properties of the uniform distribution 
imply the product formula
for the probability that  the intervals $]P_{j},P_{j-1}[$ contain $n_j$ uniform points,
\begin{equation}\label{for-lim}
{n\choose n_1,\ldots,n_k}
p(n_1+\ldots+n_k:n_k)p(n_1+\ldots +n_{k-1}:n_{k-1})\ldots
p(n_1:n_1),
\end{equation}
where the multinomial coefficient can be factored as
$\prod_{j=1}^k
{n_1+\ldots +n_j \choose n_j}.$ 
While this formula implies in an easy way the joint distribution of the occupancy counts read right-to-left,
there is no simple formula for the joint distribution of the counts read left-to-right. We will see that
in the $n\to\infty$ limit there is a considerable simplification, as in \cite{selfsim}.

\par Observe that $Z_n:=(Z_n^{(i)}: i=1,2,\ldots)$ can be
defined in the same `balls-in-boxes' fashion in terms of the
inflated  sets $nB$ and ${\cal U}_n:=\{nU_{n,j}:~1\leq j\leq n\}$.
From the extreme-value theory we know that, as $n\to\infty$, the
point process ${\cal U}_n$  converges vaguely to a unit Poisson
process $\cal U$ on $\mathbb{R}_+$, Here and henceforth, the vague
convergence means weak convergence on every finite interval
bounded away from $0$. On the other hand, $nB$ also converges
vaguely to some point process $\cal B$ on ${\mathbb R}_+$ which is
self-similar, that is satisfies $c{\cal B}=_d {\cal B}$ for every
$c>0$. The convergence of $nB$ is a consequence of the classical
renewal theorem applied to the finite-mean random walk $\{-\log
P_j: j\in\mn_0\}$. The self-similarity in this context is
analogous to the stationarity in the (additive) renewal theory.

\par The set ${\mathbb R}_+\setminus{\cal B}$ is itself a collection of open intervals (`boxes')
which accumulate in some way the points of $\cal U$ (`balls'),
hence we can define a nonnegative sequence of counts of
`balls-in-boxes' $Z:=(Z^{(i)}: i=1,2,\ldots)$ which starts with
some positive number $Z^{(1)}$ of Poisson points falling in the
leftmost nonempty interval. In view of the convergence of the
point processes, one can expect that the convergence of the
counting sequences also holds.

\begin{thm}\label{main1}
 As $n\to\infty$,
\begin{equation}\label{wc}
(Z_n^{(1)}, Z_n^{(2)},
\ldots)\to_d (Z^{(1)}, Z^{(2)},\ldots).
\end{equation}
The distribution of the limit sequence is given by the formula
\begin{eqnarray}\label{defin}
\mmp\{Z^{(1)}=n_1, \ldots, Z^{(\ell)}=n_\ell\}=\notag\\
{1\over \mu(n_1+\ldots+n_\ell)}{n_1+\ldots+n_\ell\choose n_1,\ldots,n_\ell}
p(n_1+\ldots+n_\ell:n_\ell)p(n_1+\ldots+n_{\ell-1}:n_{\ell-1})\ldots p(n_1:n_1)
\end{eqnarray}
for any  $\ell>0$, $n_1>0$ and $n_2,\ldots,n_\ell\geq 0$.
\end{thm}

\begin{proof}
Fix $\epsilon>0$ and restrict all point processes to $[\epsilon, \epsilon^{-1}]$.
By Skorohod's theorem we can select probability space in such a way that the convergence
of $(nB,{\cal U}_n)$ to $({\cal B},{\cal U})$
holds almost surely, then for  continuity reasons (see Lemma \ref{conv} to follow) 
the occupancy numbers of the intervals within $[\epsilon, \epsilon^{-1}]$ converge.
The weak convergence (\ref{wc})
follows by sending $\epsilon\to 0$ and noting that the probability that
any $m$ leftmost points of $\cal U$ fit in $[\epsilon, \epsilon^{-1}]$ goes to one.

\par Let $n=n_1+\ldots+n_\ell$ and denote by $X$ the $(n+1)$st leftmost point of $\cal U$.
The generic sequence of occupancy numbers which gives rise to the
event in (\ref{defin}) is of the form
$(n_1,\ldots,n_\ell,0,\ldots,0,m)$ where $m$ is some positive
number and the number of $0$'s is arbitrary. Let $G=\max({\cal
B}\cap [0,X])$ be the largest point of $\cal B$ smaller than $X$;
from selfsimilarity and \cite{selfsim} we know that the
distribution of $G/X$ has density $(\mu x)^{-1}{\mathbb P}\{W<x\}$
on $[0,1]$, and from the order statistics property of the Poisson
process we know that given $X$ the first $n$ points of $\cal U$
are distributed as a uniform sample from $[0,X]$. The pattern
$(n_1,\ldots,n_\ell,0,\ldots,0,m)$ occurs when the uniform
$n$-sample does not hit $[G,X]$ (event $E_1$) and within $[0,G]$
the occupancy numbers are $(n_1,\ldots,n_\ell,0,\ldots,0)$ (event
$E_2$). Integrating by parts, the probability of $E_1$ is
$$\int_0^1 {x^n{\mathbb P}\{W<x\}\over\mu x}\,{\rm d}x={1\over \mu n} \left(1-\me[W^n]\right).$$
For $i=0,1,\ldots$ let $E_{2,i}$ be the event that the pattern
$(n_1,\ldots, n_l, 0,\ldots, 0)$ with exactly $i$ zeroes occurs. In
view of the equality $\left({1 \over G}   \mathcal{B}\right)\bigcap [0,1]=_d B\setminus\{1\}$,
the conditional probability $\mmp\{E_{2,i}|X=x, E_1 \}$ equals the
probability (\ref{for-lim}) with $k=\ell+i$ and $n_{\ell+1}=\cdots=n_k=0$.
Since $E_2=\bigcup_{i=0}^\infty E_{2,i}\,$,
summing the last probabilities over $i$, we have
$$\mmp\{E_2|X=x, E_1\}={1\over 1- \me[W^n]}{n\choose n_1,\ldots,n_\ell}
p(n_1+\ldots+n_\ell:n_\ell)p(n_1+\ldots+n_{\ell-1}:n_{\ell-1})\ldots
p(n_1:n_1)=$$$$=\mmp\{E_2|E_1\}.$$ Since the probability in
(\ref{defin}) equals $\mmp\{E_1\bigcap E_2\}$, the proof is
complete.
\end{proof}

\section{$r$-counts}

We wish to connect the asymptotics of $r$-counts to Theorem
\ref{main1}. Let $Y$ be the leftmost atom of $\cal U$. For $r\geq 0$
let $K_r^*$ be the number of intervals of $]Y,\infty[ \,\setminus
{\cal B}$ that contain exactly $r$ points of $\cal U$. For $r>0$
we can take $]0,\infty[$ instead of $]Y,\infty[$ in this
definition.

\begin{lemma} \label{conv}
Let $A,B$ be two simple (i.e. without multiple points)
point pro\-ces\-ses defined and a.s. finite in some interval
$[s,t]$, and such that $A\cap B=\varnothing$ a.s. Suppose we have
weak convergence $(A_n,B_n)\to_d (A,B)$ for a sequence of
bivariate point processes. Define a gap to be a subinterval of
$[s,t]$ whose endpoints are consecutive atoms of $B$. Let $L_k$ be
the number of gaps in $B$ that contain exactly $k$ points of $A$
(with the convention that $L_0$ counts the gaps to the right of
the leftmost $A$-point  in $[s,t]$), and let $L_{n, k}$ be defined
similarly in terms of $(A_n,B_n)$. Then
$(L_{n,0},L_{n,1},\ldots)\to_d (L_{0},L_{1},\ldots)$ as
$m\to\infty$.
\end{lemma}
\begin{proof}
By Skorohod's theorem a version of the processes can be
defined on some probability space in such a way that with
probability one the convergence is pointwise. That is to say, for
large enough $n$, $\#B_n$ and $\#B$ are equal and the
points of $B_n$ (labelled, e.g. in the increasing order) are
$\epsilon$-close to the points of $B$. Same for $A_n, A$. Thus for
large $n$, there is a bijection between the gaps in $B$ and in
$B_n$ and between the points of $A$ and $A_n$ that fall in each
particular  gap.
\end{proof}
\noindent A variation of the lemma allows accumulation of atoms of
the gaps-generating process at the left endpoint of the underlying
interval. In our situation  both $\cal B$ and $\cal U$ live on the
halfline and accumulate at infinity, hence to pass from the
occupancy counts to $K_{n,r}$'s we need to take further care by
showing that the contribution  of the counts within $[s,\infty]$
is for large $s$  negligible. To this end, it is enough to work
with expected values.

\par Now, the mean contribution of $[0,s]$ to $\me[K_r^*]$ can be estimated by the expected number of points
in ${\cal B}\cap [\min(Y,s),s]$,
$$\me\left[\int_{\min(Y,s)}^\infty {{\rm d}x\over \mu x}\right]=\int_0^\infty e^{-z}{\rm d}z\int_{\min(z,s)}^s {{\rm d}x\over \mu x}
=\int_0^s e^{-z}{\rm d}z\int_z^s {{\rm d}x\over \mu
x}<\infty.$$

\begin{lemma}\label{exp}
 We have $\me[K_r^*]=(\mu r)^{-1}$  for $r>0$, 
and also $\me[K_0^*]=\nu/\mu$, where
$$\nu:=\me[-\log (1-W)]$$
may be finite or infinite.
\end{lemma}

\begin{proof}

Indeed, by the renewal theory the intensity measure of the process $\cal B$ is $(\mu x)^{-1}{\rm d}x$. Understanding
a possible $\cal B$-atom in ${\rm d}x$ as the {\it right} endpoint of a gap
we obtain for $r>0$
\begin{eqnarray*}
\me[K_r^*] =\me\left[\int_0^\infty   e^{-x(1-W)}~{x^r(1-W)^r\over r!} {{\rm d}x\over \mu x}   \right]=
{1\over \mu\, r!}\,\me\left[ \int_{0}^\infty e^{-y}y^{r-1}{\rm d}y
\right]={1\over \mu r}\,.
\end{eqnarray*}
For  $r>0$ setting $s=0$ we obtain  $\me[K_r^*]=(\mu r)^{-1}$.
For $r=0$  we have
$${\mathbb E}[K^*_0]=\me\left[\int_{0}^\infty {e^{-(1-W)x}(1-e^{-Wx})}{{\rm d}x\over \mu x}\right]={{\mathbb E}[-\log(1-W)]\over \mu}={\nu\over \mu}\,,$$
where the second factor in the integrand stands for the event that $X$ is smaller than the left endpoint of the gap.
\end{proof}

\noindent In the case $\nu=\infty$ the source of divergence of
$K_{0}^*$ is $\infty$ and not $0$, as one sees by checking that
the mean number of $\cal B$-points in $[\min(Y,s),s]$ is finite
for every $s>0$.

\begin{assertion}
The conditions $\nu=\infty$ and $K_0^*=\infty$ a.s. are equivalent.
\end{assertion}
\begin{proof}
If $K_{0}^*=\infty$ a.s. then 
$\nu=\infty$ by Lemma
\ref{exp}. 
The proof in the other direction follows by application of the Kochen-Stone extension of the Borel-Cantelli lemma.
\end{proof}

\begin{thm} As $n\to\infty$ we have
$$(K_{n,0},K_{n,1},\ldots)\to_d (K_{0}^*,K_{1}^*,\ldots),$$
along with the convergence of expectations
$${\mathbb E}[K_{n,r}]\to {\mathbb E} [K_r^*],$$
where the limit may be finite or infinite for $r=0$.
\end{thm}
\begin{proof}
The limit set satisfies  ${\cal B}\cap [0,1]=_d W_0 B$ where $B$ and $W_0$ are independent,
and $W_0$ has the density $(\mu x)^{-1}{\mathbb P}(W<x){\rm d}x$ on $[0,1]$.
We shall speak of  $[r]$-counts meaning
the intervals within $[U_{n,1},1]$ (or $[X,\infty]$, depending on the context)
that contain at most $r$ sampling points,
and we denote
$K_{[r]}^*=\sum_{i=0}^r K_r^* ,~ K_{n,[r]}=\sum_{i=0}^r K_{n,r}$.
Replacing  in the proof of Lemma \ref{exp} the lower limit of integration  $0$ by $s W_0$ we see that choosing
$s$ large enough we can achieve that the contribution to ${\mathbb E}[K_{[r]}^*]$ of the intervals with right endpoint in  $[s W_0,\infty]$
is arbitrarily small.
It remains to show that the contribution to ${\mathbb E}[K_{n,[r]}]$ of $[s/n,1]$ is small for large enough $s$ uniformly in $n$.

\par Observe that the number of components of $B^c\cap [\epsilon,1]$ that contain no more than $r$ uniform points
is nonincreasing with $n$, because the number of `balls' in a `box' can only grow as more `balls' are thrown.
Furthermore, observe that, for the purpose of estimate,
the fixed-$n$ uniform sample can be replaced by the Poisson sample of rate $n$ on $[0,1]$.
Indeed,
the probability that a gap of size $x$ is hit by $r$ uniform points  is ${n\choose r}x^r(1-x)^{n-1}$,
and in the possonised model it is $e^{-nx} (nr)^r/r!$  for $r>0$, while
 for $r=0$ we have $(1-x)^n$ versus $e^{-nx}$.
In the range $1/2<x<1$ we have elementary estimates  $c_1e^{-x}< 1-x < c_2 e^{-x}$ for suitable positive $c_1,c_2$,
which  allow to show that the mean number of $[r]$-counts coming from  $[s/n,1]$ is of the same order for both models.
The intervals of size larger $1/2$ can be ignored, since the probability that they accomodate $r$ or less sample points decays exponentially  with $n$.

\par Arguing within the framework of Poisson sample ${\cal U}\cap [0,1]$,
we compare occupancy of `boxes' generated by $B$ with that for $W_0B$. The `meander interval' $[W_0,1]$ gives negligible
contribution to $[r]$-counts hence will be ignored.
Because ${\cal B}\cap [W_0 s/n, W_0]$ is a zoomed-in copy of $[s/n,1]$,
the  sequence of occupancy counts for  ${\cal B}\cap [W_0 s/n, W_0]$ has
the same distribution  
as if we had $B\cap [s/n, 1]$ in the role of `boxes' and
a mixed Poisson process with rate $nW_0$ in the role of `balls'. By monotonicity and because $W_0<1$, the number of $[r]$-counts derived from ${\cal B}\cap [W_0s/n, W_0]$ is larger
than the number of $r$-counts  from $B\cap [s/n,W_0]$, therefore  the mean number of such counts
can be kept small by the choice of $s$. This implies the desired estimate of
the contribution of $[s/n,W_0]$ to ${\mathbb E}K_{n,[r]}$.
\end{proof}



\end{document}